\documentclass[leqno]{article}
\usepackage{amssymb, amsmath, amsfonts}

\newtheorem{theorem}{Theorem}[section]
\newtheorem{corollary}[theorem]{Corollary}
\newtheorem{lemma}[theorem]{Lemma}
\newtheorem{proposition}[theorem]{Proposition}

\numberwithin{equation}{section}

\def\sqr#1#2{{\vcenter{\vbox{\hrule height.#2pt
    \hbox{\vrule width.#2pt height#1pt \kern#1pt
    \vrule width.#2pt}
    \hrule height.#2pt}}}}

\def\ga{\gamma}
\def\al{\alpha}
\def\De{\Delta}

\def\la{\lambda}
\def\eps{\varepsilon}
\def\ka{\kappa}

\def\om{\omega}
\def\Om{\Omega}

\def\pa{\partial}

\def\bB{\bar{B}}

\def\bv{\bar{v}}

\font\bbb=msbm10

\def\R{\hbox{{\bbb R}}}

\def\capa{\, {\rm cap}\, }
\def\mes{\,{\rm mes}\, }
\def\Lip{\, {\rm Lip}}

\def\Int{{\rm Int}\,}

\def\cD{{\mathcal D}}
\def\cE{{\mathcal E}}

\def\ms{\medskip}

\newcommand{\meanint}{{\int{\mkern-18mu}-}}

\begin{document}
%\input{esint.sty}
%\fint %version of \meanint
%$\meanint f(x)dx$
%%%%%%%%%%%%%%%%%%%%%%%%%%%%%%%%%%%%%%%%%%%%%%%%%%%%%%%%%%%%%%%%%%%
\title{Can one see the fundamental frequency\\ of a drum?}

\author
{
{\bf Vladimir Maz'ya}\\
Ohio State University,
Columbus, OH, USA
\\
%Department of Mathematical Sciences, 
University of Liverpool, Liverpool,  UK
\\
%Department of Mathematics
Link\"oping University,
Sweden
%E-mail: vlmaz@math.ohio-state.edu, vlmaz@mai.liu.se
\bigskip\\
{\bf Mikhail Shubin}
\thanks{Research partially supported by NSF grant DMS-0107796} 
\\
%Department of Mathematics
Northeastern University, Boston, USA
%E-mail: shubin@neu.edu
}

\date{}

\maketitle
%%%%%%%%%%%%%%%%%%%%%%%%%%%%%%%
%                                                                 %
%       Please insert now the article body.                       %
%                                                                 %
%%%%%%%%%%%%%%%%%%%%%%%%%%%%%%%

%\subjclass{35P15}

%\keywords{Dirichlet Laplacian, bottom of spectrum, capacity}

\begin{center}
\textit{
%Dedicated 
To the memory of Felix A. Berezin 
}
\end{center}

\bigskip
\begin{abstract}
We establish two-sided estimates for the fundamental frequency
(the lowest eigenvalue) of the Laplacian in an open set $\Om\subset \R^n$ 
with the Dirichlet boundary condition. This is done
in terms of the interior capacitary radius of $\Om$ which is
defined as the maximal possible radius of  a ball $B$ with a \textit{negligible} 
intersection with the complement of $\Om$.
Here negligibility of $F\subset B$ means that
$
\capa (F)\le \ga \capa (B),
$
where $\capa$ means the Wiener (harmonic) capacity and  $\ga$ is arbitrarily fixed with the sole restriction $0<\ga <1$. 
We provide explicit values of constants in the two-sided estimates.
%For sufficiently small $\ga>0$ this result was established in  
%\cite{Mazya-74}.
%This substantially improves the result of  \cite{Mazya-74} where only sufficiently small $\ga$ were %considered.
\end{abstract}


\section{Main result}\label{S:main}
  
 Let  us consider an open set $\Omega\subset\R^n$ and denote
 the bottom of the spectrum of its minus Dirichlet Laplacian $(-\De)_{{\rm Dir}}$ by $\la(\Om)$.
 (We understand the minus Dirichlet Laplacian as the self-adjoint operator
 which is the Friedrichs extension of the operator $-\De$ defined on $C_0^\infty(\Om)$. )
 In case when $\Om$ is a bounded domain with a sufficiently regular boundary,
 $\la(\Om)$ is the lowest eigenvalue of $-\De$ with the Dirichlet
 boundary condition on $\pa\Om$. In the general case we can write
 \begin{equation}\label{E:la}
\la(\Om)= \inf_{u\in C_0^\infty(\Om)}
\frac{\int_\Omega |\nabla u|^2 dx}{\int_\Omega |u|^2 dx}.
\end{equation}
 It follows that $\Om'\subset\Om$ implies $\la(\Om)\le \la(\Om'))$.
 In particular, if $B_r$ is an open ball of radius $r$, such that $B_r\subset\Om$,
 then $\la(\Om)\le \la(B_r)=C_n r^{-2}$ where $C_n=\la(B_1)$.
 It follows that for the interior radius of $\Om$, which is defined  as
 \begin{equation*}
 r_\Om=\sup\{r|\,\exists B_r\subset\Om\},
 \end{equation*}
 we have
  \begin{equation*}
 \la(\Om)\le C_nr_\Om^{-2}.
  \end{equation*}
 But this estimate is not good for unbounded domains or
 domains with complicated boundaries.
 In particular, a similar estimate from below does not hold
 in general.
 
 The way to improve this estimate is to relax the requirement
 for $B_r$ to be completely inside $\Om$ by allowing
 some part of $B_r$, which has a ``small" Wiener capacity, to stick
 out of $\Om$. Namely, let us take an {\it arbitrary} $\ga\in (0,1)$
 and call a compact set $F\subset \bB_r$ negligible
 (or, more precisely, $\ga$-negligible) if
 \begin{equation}\label{E:ga-negl}
 \capa(F)\le\ga\capa(\bB_r).
 \end{equation}
 (Here $\capa(F)$ denotes the Wiener capacity of $F$, $\bB_r$ is the closure of $B_r$.)
 
 Now denote
 \begin{equation*}
 r_{\Om,\ga}=\sup\{r|\, \exists B_r, \ \bB_r\setminus\Om \textrm{ is } 
 \ga\textrm{-negligible}\}.
 \end{equation*}
This is the \textit{interior capacitary radius}.

 \begin{theorem}\label{T:main} Let us fix $\ga\in(0,1)$. Then 
 there exist $c=c(\ga, n)>0$ and $C=C(\ga, n)>0$, such that for every 
 open set $\Om\subset \R^n$
 \begin{equation}\label{E:main}
 cr_{\Om,\ga}^{-2}\le \la(\Om)\le Cr_{\Om,\ga}^{-2}.
 \end{equation}
 Explicit values of constants $c=c(\ga,n)$ and $C=C(\ga, n)$
 are provided below in \eqref{E:c-explicit} and \eqref{E:C-ga-n}
 respectively. 
 \end{theorem}
 Let us formulate some  corollaries of  this theorem.

\begin{corollary}\label{C:1}
$\lambda(\Omega)>0 \textrm{ if and only if } r_{\Omega,\gamma}<\infty$.
\end{corollary}

This corollary gives a necessary and sufficient condition of strict positivity
of the operator $(-\De)_{\rm Dir}$ in $\Om$.

Since the condition $\lambda(\Omega)>0$ does not contain $\ga$,
we immediately obtain
\begin{corollary}\label{C:2}
Conditions $r_{\Omega,\gamma}<\infty$, taken for different
$\gamma$'s, are equivalent.
\end{corollary}

Denoting $F=\R^n\setminus \Om$ (which can be an arbitrary closed
subset in $\R^n$), we obtain from the previous corollary
(comparing $\ga=0.01$ and $\ga=0.99$):

\begin{corollary} \label{C:3} 
Let $F$ be a closed subset in $\R^n$, which has
the following property: there exists $r>0$ such that
$$
\capa(F\cap\bB_r)\geq 0.01\capa(\bB_r)
$$
for all $B_r$.  Then there exists $r_1>0$ such that
$$
\capa(F\cap\bB_{r_1})\geq 0.99\capa(\bB_{r_1})
$$
for all $B_{r_1}$.
\end{corollary}

This is  a new property of capacity which is proved by
spectral theory arguments.

Once upon a time Marc Kac \cite{Kac} formulated a  
fascinating question:
\textbf{``Can one hear the shape of a drum?"}  The precise meaning
of this question is as follows: is it possible to reconstruct the drum
(a bounded domain  in $\R^2$) up to an isometry by the spectrum 
of its  Dirichlet Laplacian?

Theorem \ref{T:main} suggests formulation of a question,
which is roughly inverse to the question above:
\textbf{``Can one see the fundamental frequency of a drum?"}
More precisely, can one find a simple visual image
related to a domain in $\R^2$ (or $\R^n$), such that it allows to recover
the lowest eigenvalue of the Dirichlet Laplacian in this domain,
or at least give reasonably good estimates of this eigenvalue?
Assuming that our eye (possibly armed by a visual aid device)
can filter out the sets of  ``small" capacity,
a partial answer to this question is given by Theorem \ref{T:main}.
 
The inequalities \eqref{E:main} for sufficiently small $\ga>0$
(comparable with $(4n)^{-4n}$) 
were established in \cite{Mazya-74} (see also Chapters 10 and 11 in \cite{Mazya8}). 
Theorem \ref{T:main} provides a substantial improvement,
in particular allowing corollaries \ref{C:2} and \ref{C:3}
and providing explicit values of the constants.

The proof of Theorem \ref{T:main} is based 
on the ideas of our paper \cite{Mazya-Shubin}.
Necessary definitions and
results about the Wiener capacity can be found
e.g. in \cite{Mazya8, Wermer}. 

\ms
\textbf{Acknowledgments.} We are grateful to Egon Schulte
who communicated to us the coverings multiplicity estimate \eqref{E:cov},
and also to Dan Grieser and Elliott Lieb for useful references.

\section{Preliminaries on capacity}\label{S:nec}
In this section we will  recall some  definitions and introduce necessary notations.
For simplicity we will always assume that $n\ge 3$.

For every subset $\cD\subset\R^n$ denote by $\Lip(\cD)$ the space of 
(real-valued) functions satisfying
the uniform Lipschitz condition in $\cD$, and by $\Lip_c(\cD)$ the subspace in $\Lip(\cD)$
of all functions with compact support in $\cD$ (this will be only used when $\cD$ is open).
By $\Lip_{loc}(\cD)$ we will denote the set of functions on (an open set) $\cD$ which are
Lipschitz on any compact subset $K\subset\cD$.
Note that $\Lip(\cD)=\Lip(\bar\cD)$ for any bounded $\cD$.

If $F$ is a compact subset in  ${\R}^n$, then the 
Wiener capacity of $F$ 
is defined as
\begin{equation}\label{E:cap-def}
{\capa}(F)=\inf\left\{\left.
\int_{\R^n}|\nabla u(x)|^2 dx\,\right|\; u\in\Lip_c(\R^n), u|_F=1\right\}.
\end{equation}

Note that the infimum does not change if we  restrict ourselves to the  
functions $u\in\Lip_c(\R^n)$ such that $0\le u\le 1$  everywhere 
(see e.g. \cite{Mazya8}, Sect. 2.2.1).

\ms
We will also need another (equivalent) definition of the Wiener capacity $\capa(F)$ 
for a compact set $F\subset \bB_r$. For  $n\ge 3$ it is as follows:
\begin{equation}\label{E:cap-def2}
\capa(F)=\sup\{\mu(F)\left| \int_F\cE(x-y) d\mu(y)\le 1\quad \text{on}\ \R^n\setminus F\right.\},
\end{equation}
where the supremum is taken over all positive finite Radon measures $\mu$ on $F$ and 
$\cE=\cE_n$ is the standard fundamental solution of $-\De$ in $\R^n$ i.e.
\begin{equation*}\label{E:fund-sol}
\cE(x)=\frac{1}{(n-2)\om_n}|x|^{2-n}\;,
\end{equation*}
with $\om_n$ being the area of the unit sphere $S^{n-1}\subset\R^n$. 
The maximizing measure in \eqref{E:cap-def2} exists and is unique.
We will denote it $\mu_F$ and call it the \textit{equilibrium measure}.
Note that
\begin{equation}\label{E:mu-F}
\capa(F)
= \mu_F(F)
=\mu_F(\R^n)
=\langle\mu_F,1\rangle
=\int_F d\mu_F.
\end{equation}
The corresponding potential will be denoted $P_F$, so 
\begin{equation*}\label{E:psi-F}
P_F(x)=\int_F\cE(x-y)d\mu_F(y), \quad x\in \R^n\setminus F.
\end{equation*} 
We will call $P_F$ the \textit{equilibrium potential} or \textit{capacitary potential}. 
We will  extend it to $F$ by setting $P_F(x)=1$ for all $x\in F$.

It follows from the maximum principle that $0\le P_F\le 1$
everywhere in $\R^n$.

In case when  $F$ is the closure
of an open subset with a smooth boundary, $u=P_F$ is the unique minimizer
for the Dirichlet integral in \eqref{E:cap-def}.
 In particular, 
\begin{equation*}\label{E:minimizer}
\int |\nabla P_F|^2 dx=\capa(F).
\end{equation*}
where the integration is taken over $\R^n$ (or $\R^n\setminus F$).

The capacity of the ball $\bB_r$ is easily calculated and is given by
\begin{equation}\label{E:cap-ball}
\capa(\bB_r)=(n-2)\om_n r^{n-2}.
\end{equation}

\section{Lower bound}\label{S:lower}

In this section we will establish the lower bound for $\la(\Om)$ from Theorem \ref{T:main} which 
 is an easier part of this theorem.  The key part of the lower bound proof  is presented in the following lemma,  
which was first proved in \cite{Mazya-63b} (see also  \cite{Mazya8},  where it is 
present  as  a particular case of
a much more general Theorem 10.1.2, part 1),  though without an explicit constant,
which we provide to specify explicit constants in Theorem \ref{T:main}.

\begin{lemma}\label{L:Mazya1}
The following inequality holds for every complex-valued
function $u\in \Lip(\bB_r)$  which
vanishes on a compact set $F\subset \bB_r$ (but is not identically $0$ on $\bB_r$):
\begin{equation}\label{E:cap-above}
\capa(F)\le \frac{C_n\int_{B_r} |\nabla u(x)|^2dx}{r^{-n}\int_{B_r} |u(x)|^2 dx}\;,
\end{equation}
where
\begin{equation}\label{E:Mazya1-Cn}
C_n = 4\om_n\left(1-\frac{2}{n^2}\right)
\end{equation}
\end{lemma}
\textbf{Beginning of Proof.} 
\textbf{A.} Clearly, it is sufficient to consider the ball $B_r$ centered at $0$,
and real-valued functions $u\in \Lip(\bB_r)$. By scaling we see that it suffices to consider
the case $r=1$. (The corresponding estimate for an arbitrary $r>0$ follows from the one 
with $r=1$ with the same constant $C_n$.) So we need to prove the estimate
\begin{equation}\label{E:cap-above1}
\int_{B_1}|u|^2 dx \le
\frac{C_n}{\capa(F)}\int_{B_1}|\nabla u|^2 dx,
\end{equation}
where $F$ is a compact subset of $\bB_1$, $u\in\Lip(\bB_1)$ and $u|_F=0$.

To be able to use \eqref{E:cap-def}, consider the following function $U\in\Lip(\R^n)$:
\begin{equation*}\label{E:U}
U(x)=
\begin{cases}
1-|u(x)|,  
& \text{if $|x|\le 1$,}\\
|x|^{2-n}(1-|u(|x|^{-2}x)|),
& \text{if $|x|\ge 1$},
\end{cases}
\end{equation*}
i.e. $U$ extends $1-|u|$ to $\{x:|x|\ge 1\}$ as the Kelvin transform of $1-|u|$.
Clearly, $U|_F=1$, $|\nabla U|=|\nabla u|$ almost everywhere in $B_1$, 
$U(x)=O(|x|^{2-n})$ and $|\nabla U(x)|=O(|x|^{1-n})$ as $|x|\to\infty$. It follows that
$U$ can serve as a test function in \eqref{E:cap-def}, i.e.
\begin{equation}\label{E:cap-U}
\capa(F)\le \int_{\R^n}|\nabla U|^2 dx.
 \end{equation}
Using the harmonicity of $|x|^{2-n}$ and the Green-Stokes formula, we obtain
by a straightforward calculation
\begin{equation}\label{E:nabla-U}
\int_{\R^n}|\nabla U|^2 dx
= 2\int_{B_1}|\nabla u|^2 dx + (n-2) \int_{\pa B_1}(1-|u(\om)|)^2 d\om,
\end{equation}
where $\pa B_1=\{\om\in\R^n,\ |\om|=1\}$ is the boundary of $B_1$
(the unit sphere in $\R^n$), $d\om$ means the standard  volume
element on $\pa B_1$.

\bigskip
\noindent
\textbf{B.}
For a function $v$ on $\pa B_1$ define its average 
\begin{equation*}
\bar v=\meanint_{\pa B_1}v d\om
=\frac{1}{\om_n}\int_{\pa B_1}v d\om\,.
 \end{equation*}

To continue the proof of Lemma \ref{L:Mazya1},
 we will need the following elementary lemma.
 \begin{lemma}\label{L:v-bv} 
 For any $v\in\Lip(B_1)$,
 \begin{equation}\label{E:v-bv}
 \int_{\pa B_1}|v-\bv|^2 d\om
\le \int_{B_1}|\nabla v|^2 dx.
 \end{equation}
\end{lemma}

\textbf{Proof of Lemma \ref{L:v-bv}.} 
It suffices to prove it for real-valued functions $v$. 
Let us expand $v$ in spherical functions. Let 
\begin{equation*}\label{E:Ykl}
\{Y_{k,l}|\, l=0,1, \dots, n_k, \ k=0,1,\dots\}
\end{equation*}
be an orthonormal basis in $L^2(\pa B_1)$ which consists 
of eigenfunctions of the (negative) Laplace-Beltrami
operator $\De_\om$ on $\pa B_1$, so that the eigenfunctions $Y_{k,l}=Y_{k,l}(\om)$
with a fixed $k$ have the same eigenvalue $-k(k+n-2)$ (which has multiplicity $n_k+1$).
Note that the zero eigenvalue (corresponding to $k=0$) has multiplicity 1 and 
$Y_{0,0}=const=\om_n^{-1/2}$ for the corresponding eigenfunction.

Writing $x=r\om$, where $r=|x|$, $\om=x/|x|$, we can present $v$ in the form
\begin{equation}\label{E:v-expansion}
v(x)=v(r,\om)=\sum_{k,l}v_{k,l}(r)Y_{k,l}(\om).
\end{equation}
Then
\begin{equation}\label{E:L2-norm}
\int_{B_1}|v(x)|^2 dx = 
\sum_{k,l}\int_0^1 |v_{k,l}(r)|^2  r^{n-1} dr,
\end{equation}
and
\begin{equation}\label{E:L2-norm-s}
\int_{\pa B_1}|v(\om)|^2 d\om = 
\sum_{k,l} |v_{k,l}(1)|^2.
\end{equation}
It follows that
\begin{equation}\label{E:v-bv2}
\int_{\pa B_1}|v(\om)-\bv|^2 d\om = 
\sum_{\{k,l: k\ge 1\}}|v_{k,l}(1)|^2.
\end{equation}
Taking into account that
\begin{equation*}
|\nabla v|^2 = \left|\frac{\pa v}{\pa r}\right|^2 + r^{-2} |\nabla_\om v|^2,
\end{equation*}
where 
$\nabla_\om$ means the gradient along the unit sphere
with variable $\om$ and fixed $r$, we also get
\begin{equation}\label{E:Dir-int}
\int_{B_1}|\nabla v|^2 dx =
\sum_{k,l}\int_0^1 \left(|v'_{k,l}(r)|^2 + \frac{k(k+n-2)}{r^2}|v_{k,l}(r)|^2 \right)r^{n-1} dr.
\end{equation}
Comparing \eqref{E:v-bv2} and \eqref{E:Dir-int}, and taking into account that
$k(k+n-2)$ increases with $k$, we see that
it suffices to establish that the inequality
\begin{equation*}\label{E:k-eq-1}
|g(1)|^2 \le \int_0^1 \left( |g'(r)|^2 + \frac{n-1}{r^2} |g(r)|^2 \right) r^{n-1} dr, 
\end{equation*}
holds for any real-valued function  $g\in \Lip([0,1]))$.  To this end write
\begin{align*}
g(1)^2
&=\int_0^1(r^{n-2}g^2)' dr
=\int_0^1 [2r^{n-2}g'g + (n-2)r^{n-3} g^2] dr\\
&\le \int_0^1 [r^{n-1}g'^2 +(n-1)r^{n-3}g^2] dr
=\int_0^1\left( g'^2 + \frac{n-1}{r^2}g^2\right) r^{n-1}dr,
\end{align*}
which proves Lemma \ref{L:v-bv}. $\square$

\bigskip
\noindent
\textbf{C. Proof of Lemma \ref{L:Mazya1} (continuation).}
Let us normalize $u$ by requiring $\overline{|u|}=1$,
i.e. average of $|u|$ over $\pa B_1$ equals 1.  This can be done if 
$u\not\equiv 0$ on $\pa B_1$.
Then by Lemma \ref{L:v-bv} 
\begin{equation*}
\int_{\pa B_1}(1-|u|)^2 d\om
\le \int_{B_1} |\nabla u|^2 dx. 
\end{equation*}
Combining this with \eqref{E:cap-U} and \eqref{E:nabla-U}, we obtain
\begin{equation*}
\capa(F)
\le n \int_{B_1} |\nabla u|^2 dx. 
\end{equation*}
Removing the restriction $\overline{|u|}=1$, we can conclude that
for any $u\in \Lip(B_1)$
\begin{equation}\label{E:ou2}
\left(\meanint_{\pa B_1} |u| d\om \right)^2
\le \frac{n}{\capa(F)}\int_{B_1} |\nabla u|^2 dx. 
\end{equation}
(This  obviously also holds in case when $u\equiv 0$ on $\pa B_1$.)

Note that for any real-valued function $v\in \Lip(B_1)$
\begin{equation*}
\meanint_{\pa B_1}|v-\bv|^2 d\om
=\meanint_{\pa B_1} |v|^2 d\om - \bv^2,
\end{equation*}
 hence, using \eqref{E:v-bv}, we get
 \begin{equation*}
 \meanint_{\pa B_1} |v|^2 d\om 
 = \bv^2 + \meanint_{\pa B_1} |v-\bv|^2 d\om
 \le \bv^2 + \frac{1}{\om_n}\int_{B_1}|\nabla v|^2 dx.  
 \end{equation*}
Applying this to $v=|u|$ and using \eqref{E:ou2},
we obtain
\begin{equation}\label{E:u2dom}
\int_{\pa B_1}|u|^2 d\om
\le \left(1+\frac{n\om_n}{\capa(F)}\right)\int_{B_1}|\nabla u|^2 dx.
\end{equation}

\bigskip
\noindent
\textbf{D.} Note that  out goal is an estimate which is similar to \eqref{E:u2dom} 
but with the integral over $\pa B_1$ in the left hand side replaced by the
integral over $B_1$.  To this end we will again use the expansion
\eqref{E:v-expansion} of $v=|u|$ over spherical functions, and the
identities \eqref{E:L2-norm},   \eqref{E:L2-norm-s} and  \eqref{E:Dir-int}.   
Let us take a real-valued function $g\in\Lip([0,1])$ and denote  
\begin{equation*}
Q=\int_0^1g^2(r)r^{n-1}dr.
\end{equation*}
Integrating by parts, we obtain
\begin{equation*} 
Q = -\frac{2}{n} \int_0^1 gg' r^{n} dr + \frac{1}{n} g^2(1).
\end{equation*}
Using an elementary inequality $2ab\le \eps a^2 + \eps^{-1}b^2$, 
where $a,b\in\R$, $\eps>0$, and taking into account that $r\le 1$,
we obtain 
\begin{align*}
Q 
&\le \frac{1}{n} \int_0^1 \left(\eps g^2(r) 
+ \frac{1}{\eps}g'^2(r)\right)r^{n-1}dr + \frac{1}{n}g^2(1)\\
&=\frac{\eps}{n}Q + \frac{1}{n\eps}\int_0^1g'^2(r)r^{n-1}dr + \frac{1}{n}g^2(1),
\end{align*}
hence for any $\eps\in(0,n)$
\begin{equation*}
Q\le \frac{1}{(n-\eps)\eps}\int_0^1g'^2(r)r^{n-1}dr + \frac{1}{n-\eps}g^2(1).
\end{equation*}
Taking $\eps=n/2$, we obtain
\begin{equation}\label{E:Q-le}
Q\le \frac{4}{n^2}\int_0^1g'^2(r)r^{n-1}dr + \frac{2}{n}g^2(1).
\end{equation}
Now we can argue as in the proof of Lemma \ref{L:v-bv}, expanding $v=|u|$
over spherical harmonics $Y_{k,l}$. Then the desired inequality follows
from the inequalities for the coefficients $v_{k,l}=v_{k,l}(r)$, with the strongest one
corresponding to the case $k=0$ (unlike $k=1$ in Lemma \ref{L:v-bv}). Then using
the inequality \eqref{E:Q-le} for $g=v_{0,0}$ we obtain
\begin{equation}\label{E:B1-u2}
\int_{B_1} |u|^2 dx 
\le \frac{4}{n^2}\int_{B_1} |\nabla u|^2 dx + \frac{2}{n}\int_{\pa B_1} |u|^2 d\om.
\end{equation}
Using \eqref{E:u2dom}, we deduce from \eqref{E:B1-u2}:
\begin{equation}\label{E:B1-u2-2}
\int_{B_1} |u|^2 dx 
\le \left[\frac{4}{n^2}+\frac{2}{n}\left(1+\frac{n\om_n}{\capa(F)}\right)\right]
\int_{B_1} |\nabla u|^2 dx.
\end{equation}
Taking into account  the inequality
\begin{equation*}
\capa(F)\le \capa(\bB_1) = (n-2)\om_n,
\end{equation*}
we can estimate the constant in front of the integral in the right hand side of 
\eqref{E:B1-u2-2} as follows:
\begin{equation*}
\frac{4}{n^2}+\frac{2}{n}\left(1+\frac{n\om_n}{\capa(F)}\right)
\le
\frac{4\om_n}{\capa(F)}\left(1-\frac{2}{n^2}\right),
 \end{equation*}
 which ends the proof of Lemma \ref{L:Mazya1}. $\square$

\bigskip
The lower bound in \eqref{E:main} is given by
\begin{lemma}\label{L:lower} There exists $c=c(\ga,n)>0$ such that
for all open sets $\Om\subset\R^n$
\begin{equation}\label{E:lower}
\la(\Om)\ge c r_{\Om,\ga}^{-2}.
\end{equation}
\end{lemma}

\textbf{Proof.} Let us fix $\ga\in(0,1)$ and choose any $r>r_{\Om,\ga}$. 
Then any ball $\bB_r$  has a non-negligible
intersection with $\R^n\setminus\Om$, i.e. 
\begin{equation*}
\capa(\bB_r\setminus\Om) \ge \ga\capa(\bB_r).
\end{equation*}
Since any $u\in C_0^\infty(\Om)$ vanishes on $\bB_r\setminus \Om$, it follows 
from Lemma \ref{L:Mazya1} that for any such $u$
\begin{equation*}
\int_{\bB_r}|u|^2 dx
\le 
\frac{C_n}{r^{-n}\capa{(\bB_r\setminus\Om)}} \int_{\bB_r}|\nabla u|^2 dx
\le 
\frac{C_n}{r^{-n}\ga\capa{(\bB_r)}} \int_{\bB_r}|\nabla u|^2 dx.
\end{equation*}
Taking into account that $\capa(\bB_r)=\capa(\bB_1)r^{n-2}$,
we obtain
\begin{equation*}
\int_{\bB_r}|u|^2 dx
\le 
\frac{C_n r^{2}}{\ga\capa(\bB_1)} \int_{\bB_r}|\nabla u|^2 dx.
\end{equation*}
Now let us choose a covering of $\R^n$ by balls $\bB_r=\bB_r^{(k)}$, 
$k=1,2,\dots,$ so that the multiplicity  of this covering is at most $N=N(n)$. 
For example, we can make
\begin{equation}\label{E:cov}
N(n)\le n\log n + n\log(\log n) + 5n, \quad n\ge 2,
\end{equation}
which holds also for the smallest multiplicity of coverings of $\R^n$ by translations
of any convex body  (see Theorem 3.2 in \cite{Rogers}).

Then summing up the estimates above over all balls in this covering, we 
see that 
\begin{eqnarray*}
\int_{\R^n}|u|^2 dx
\le
\sum_k\int_{\bB_r^{(k)}}|u|^2 dx
\le 
\frac{C_n r^{2}}{\ga\capa(\bB_1)}\sum_k \int_{\bB_r^{(k)}}|\nabla u|^2 dx\\
\le
\frac{C_n N r^{2}}{\ga\capa(\bB_1)}\int_{\R^n}|\nabla u|^2 dx.
\end{eqnarray*}
Recalling \eqref{E:la}, we see that 
\begin{equation*}
\la(\Om)\ge cr^{-2}
\end{equation*}
with
\begin{equation}\label{E:c-explicit}
c=c(\ga,n)=\frac{\ga\capa(\bB_1)}{C_n N}
=\frac{\ga n^2 (n-2)}{4(n^2-2)N}.
\end{equation}
Taking limit as $r\downarrow r_{\Om,\ga}$, we obtain 
\eqref{E:lower} with the same $c$.
$\square$

\section{Upper bound}\label{S:upper}

\subsection{}\label{SS:4.1}
According to \eqref{E:la}, to get an upper bound for $\la(\Om)$ it is enough to 
take any test function $u\in C_0^\infty(\Om)$ and write
 \begin{equation}\label{E:la-test}
\la(\Om)\le
\frac{\int_\Omega |\nabla u|^2 dx}{\int_\Omega |u|^2 dx}.
\end{equation}
For simplicity of notations we will write $\la$ instead of $\la(\Om)$
everywhere in this section. The inequality \eqref{E:la-test} can be rewritten
as follows:
 \begin{equation}\label{E:la-test2}
\int_\Om |u|^2 dx \le \la^{-1}
\int_\Omega |\nabla u|^2 dx.
\end{equation}
By approximation, it suffices to take $u\in\Lip_c(\Om)$ 
or even $u\in H^1_0(\Om)$, where $H^1_0(\Om)$ is the closure
of $C_0^\infty(\Om)$ in the standard Sobolev space $H^1(\Om)$ 
(which consists of all $u\in L^2(\Om)$ with the
distributional derivatives $\partial u/\partial x_j\in L^2(\Om)$, $j=1,\dots, n$).

In particular, choosing a ball $B_r$, we can take
\begin{equation}\label{E:test-req} 
u\in \Lip_c(\Om\cap B_r)=\Lip_c(\Om)\cap \Lip_c(B_r). 
\end{equation}
Let us take a compact set $F\subset \bB_{3r/2}$, such that $F$ is the closure
of an open set with a smooth boundary. (In this section we will call such sets
\textit{regular} subsets of $\bB_{3r/2}$.)  
Denote by $P_F$ its equilibrium potential (see Sect. \ref{S:nec}).
Regularity of $F$ implies that $P_F\in\Lip(\R^n)$.
By definition $P_F=1$ on $F$, so $1-P_F=0$ on $F$. Let us also assume that
\begin{equation*}\label{E:F-req}
\Int F\supset\bB_r\setminus\Om, 
\end{equation*}
where $\Int F$ means the set of all interior points of $F$
(so $\Int F$  is an open subset in $\R^n$). Then $1-P_F=0$
in a neighborhood of $\bB_r\setminus\Om$. Therefore, multiplying 
$1-P_F$ by a cut-off function $\eta\in C_0^\infty (B_r)$, we will get
a function $u=\eta(1-P_F)$, satisfying the requirenment \eqref{E:test-req},
hence fit to be a test function in \eqref{E:la-test}.

In the future we will also assume that the cut-off function $\eta\in C_0^\infty (B_r)$ 
has the following properties:
\begin{equation*}\label{E:eta-prop}
0\le \eta\le 1\ {\rm on} \ B_r,
\quad
\eta=1\  {\rm on}\ B_{(1-\kappa)r},
\quad 
|\nabla\eta|\le 2(\kappa r)^{-1}\ {\rm on}\ B_r,
\end{equation*}
where $0<\kappa<1$ and the balls 
$B_r $ and $B_{(1-\kappa)r}$ are supposed to have the same center.
Using integration by parts and  the equation $\De P_F=0$ on $B_r\setminus F$, 
we obtain for the test function  $u=\eta(1-P_F)$:
\begin{align*}
\int_{B_r}|\nabla u|^2 dx  
&=\int_{B_r}\left(|\nabla\eta|^2(1-P_F)^2
-\nabla(\eta^2)\cdot (1-P_F)\nabla P_F
+\eta^2|\nabla P_F|^2  \right) dx\\
&=\int_{B_r}|\nabla\eta|^2(1-P_F)^2 dx 
\le 4(\kappa r)^{-2}\int_{B_r} (1-P_F)^2 dx.
\end{align*}
Therefore, from \eqref{E:la-test2} we obtain
\begin{equation*}
\int_{B_r}|u|^2 dx
\le \la^{-1}4(\kappa r)^{-2}\int_{B_r} (1-P_F)^2 dx.
\end{equation*} 
Since $0\le P_F\le 1$, the last integral in the right hand side is estimated by
\begin{equation*}
\mes(B_r)=n^{-1}\om_n r^n.
\end{equation*}
where $\mes$ means the usual Lebesgue measure on $\R^n$.
Therefore,
\begin{equation*}
\int_{B_r}|u|^2 dx
\le 
4 n^{-1} \om_n \la^{-1} \kappa^{-2}  r^{n-2}.
\end{equation*} 
Restricting the integral in the left hand side to $B_{(1-\ka)r}$, we obtain
\begin{equation}\label{E:1-PF}
\int_{B_{(1-\ka)r}}(1-P_F)^2 dx
\le
4 n^{-1} \om_n \la^{-1} \kappa^{-2}  r^{n-2}. 
\end{equation}

\subsection{}\label{SS:4.2}

Now we need to provide an appropriate  lower bound for the left hand side
of \eqref{E:1-PF}.  To this end we first 
restrict the integration to the spherical layer
\begin{equation*}
S_{r_1,r_2}=B_{r_2}\setminus B_{r_1}, 
\end{equation*}
where $0<r_1<r_2<r$. In the future we will take 
\begin{equation}\label{E:r1r2}
r_1=(1-2\ka)r, \quad  r_2=(1-\ka)r,
\end{equation} 
where $0<\ka<1/2$, though it is convenient to write some formulas in a bigger generality. 
Let us denote the volume of the  layer $S_{r_1,r_2}$ by $|S_{r_1,r_2}|$, i.e.
\begin{equation*}\label{E:mes-Srka}
|S_{r_1,r_2}|=\mes S_{r_1,r_2} =  n^{-1} \om_n (r_2^n - r_1^n).
\end{equation*}
We will also need the  notation 
\begin{equation*}
\meanint_{S_{r_1,r_2}}f(x)dx
=
\frac{1}{|S_{r_1,r_2}|}\int_{S_{r_1,r_2}} f(x) dx
\end{equation*}
for the average of a positive function $f$ over $S_{r_1,r_2}$.
In particular, restricting the integration in \eqref{E:1-PF} to $S_{r_1,r_2}$
(with $r_1, r_2$ as in \eqref{E:r1r2})
and dividing by $|S_{r_1,r_2}|$, we obtain
\begin{equation*}
\meanint_{S_{r_1,r_2}}(1-P_F)^2 dx
\le
\frac{4\la^{-1}\ka^{-2}r^{n-2}}{r_2^n - r_1^n}.
\end{equation*}
Hence, by the Cauchy-Schwarz inequality, 
\begin{equation}\label{E:int-square}
\left[1-\meanint_{S_{r_1,r_2}} P_Fdx\right]^2
=
\left[\meanint_{S_{r_1,r_2}} (1-P_F)dx\right]^2
\le
\frac{4\la^{-1}\ka^{-2}r^{n-2}}{r_2^n - r_1^n}.
\end{equation}
To simplify the right hand side, let us estimate $(r_2^n - r_1^n)^{-1}$ from above.
Applying the Bernoulli inequality, we see that
\begin{align*}
r_2^n - r_1^n 
= (r_2 - r_1)(r_2^{n-1}+r_2^{n-2}r_1+\dots + r_1^{n-1})\\
\ge  
n \ka r r_1^{n-1} 
=
n\ka r^n (1-2\ka)^{n-1}
\ge
n\ka r^n [1 - 2(n-1)\ka].
\end{align*}
Now note that
\begin{equation*}
\frac{1}{1-2(n-1)\ka}
\le
1+4(n-1)\ka, 
\end{equation*}
provided
\begin{equation}\label{E:ka-cond}
0<\ka \le \frac{1}{4(n-1)}.
\end{equation}
Under this condition it follows that
\begin{equation}\label{E:r2-r1}
\frac{1}{r_2^n - r_1^n}
\le
n^{-1}\ka^{-1}r^{-n} \left[1+4(n-1)\ka\right],
\end{equation}
and  \eqref{E:int-square}
takes the form
\begin{equation}\label{E:int-square2}
\left[1-\meanint_{S_{r_1,r_2}} P_Fdx\right]^2
\le
4 n^{-1} \ka^{-3} \left[1+4(n-1)\ka\right] \la^{-1} r^{-2}.
\end{equation}

\subsection{}\label{SS:4.3}
For simplicity of notations and without loss
of generality we may assume that the ball $B_r$ is centered at $0\in\R^n$
(and so are smaller balls and spherical layers).

To provide a lower bound for the left hand side
of \eqref{E:int-square2}, we
will give an upper bound for the average of
$P_F$.  According to the definition of $P_F$
and notations from Section \ref{S:nec}, we can write
\begin{align}\label{E:mean-PF-est}
\meanint_{S_{r_1,r_2}} P_F dx 
&=
\meanint_{S_{r_1,r_2}}\left( \int_F \cE(x-y) d\mu_F(y) \right) dx \\
&=
\int_F\left( \meanint_{S_{r_1,r_2}} \cE(x-y)dx \right) d\mu_F(y).
\notag
\end{align}
The inside integral in the right hand side can be explicitly calculated 
(due to Newton) as the potential of a uniformly charged spherical layer
with total charge 1.  
The result of this calculation is $|S_{r_1,r_2}|^{-1}V_{r_1,r_2}(y)$, where  
\begin{eqnarray}\label{E:layer-pot}
\qquad V_{r_1, r_2}(y)=
\begin{cases}
\frac{r_2^2-r_1^2}{2(n-2)}\, ,
&\text{if $|y|\le r_1$},
\ms
\\
- \frac{|y|^2}{2n} + \frac{r_2^2}{2(n-2)} - \frac{r_1^n}{n(n-2)|y|^{n-2}}\,, 
&\text{if $r_1\le |y| \le r_2$},
\ms
\\
\frac{r_2^n - r_1^n}{n(n-2)|y|^{n-2}}\,, 
&\text{if $|y|\ge r_2$}.
\end{cases}
\end{eqnarray}
The function $y\mapsto V_{r_1,r_2}(y)$ belongs to $C^1(\R^n)$ and is
spherically symmetric; it tends to $0$ as $|y|\to\infty$;
it is harmonic in $\R^n\setminus S_{r_1,r_2}$
and satisfies the equation $\De V_{r_1,r_2}=-1$ in $S_{r_1, r_2}$.
These properties uniquely define the function $V_{r_1,r_2}$.
Differentiating it with respect to $|y|$, we easily see that it is decreasing
with respect to $|y|$, hence its  maximum is at $y=0$
(hence given by the first row in \eqref{E:layer-pot}). So we obtain, using
\eqref{E:r2-r1}:
\begin{align*}\label{E:mean-E}
&\meanint_{S_{r_1,r_2}} \cE(x-y)dx
\le
|S_{r_1,r_2}|^{-1}V_{r_1,r_2}(0)
=
\frac{n(r_2^2 - r_1^2)}{2(n-2)\om_n(r_2^n - r_1^n)}\\
&= 
\frac{n\ka r(r_1+r_2)}{2(n-2)\om_n(r_2^n - r_1^n)}
\le
\frac{n r^2 \ka(1-\ka)}{(n-2)\om_n(r_2^n - r_1^n)}
\\
& \le 
\frac{(1-\ka)[1+4(n-1)\ka]}{(n-2)\om_n r^{n-2}}
\le
\frac{1+(4n-5)\ka}{(n-2)\om_n r^{n-2}}
\end{align*}
Finally, using \eqref{E:cap-ball}, we obtain
\begin{equation}\label{E:mean-E2}
\meanint_{S_{r_1,r_2}} \cE(x-y)dx
\le
\frac{1+(4n-5)\ka}{\capa(\bB_r)}\;.
\end{equation}
provided  $r_1, r_2$ choosen as in \eqref{E:r1r2}
and  \eqref{E:ka-cond}
is satisfied. 

\subsection{}\label{SS:4.4}

Using \eqref{E:mean-E2} in \eqref{E:mean-PF-est} and taking into account \eqref{E:mu-F},
we obtain
\begin{align}\label{E:meanint-PF-est}
\meanint_{S_{r_1,r_2}} P_F(x) dx
\le
\frac{1+(4n-5)\ka}{\capa(\bB_r)} \int_F d\mu_F(y)\\
= 
[(1+(4n-5)\ka]\frac{\capa(F)}{\capa(\bB_r)}
\le
[(1+(4n-5)\ka]\ga,
\notag
\end{align}
provided $F$ is $\ga$-negligible.
(i.e. satisfies  \eqref{E:ga-negl}).
Note  that we do not assume that $F\subset \bB_r$ but
do assume that $0<\ga<1$. In this case,
taking into account \eqref{E:ka-cond},  we can take
\begin{equation}\label{E:ka}
\ka=\min \left\{ \frac{1}{4(n-1)},\; \frac{1-\ga}{2(4n-5)\ga} \right\},
\end{equation}
so that \eqref{E:ka-cond} is satisfied, and, besides,
\begin{equation*}
[(1+(4n-5)\ka]\ga\le \frac{1+\ga}{2}=1-\frac{1-\ga}{2},
\end{equation*}
 so that \eqref{E:meanint-PF-est} becomes
\begin{equation*}\label{E:meanint-PF-final}
\meanint_{S_{r_1,r_2}} P_F(x) dx\le 1-\frac{1-\ga}{2}.
\end{equation*}
Taking this into account in \eqref{E:int-square2}
and using \eqref{E:ka-cond},
we obtain
\begin{equation*}
\frac{(1-\ga)^2}{4}
\le
4 n^{-1}\ka^{-3}[1+4(n-1)\ka] \la^{-1} r^{-2}
\le 
8 n^{-1} \ka^{-3} \la^{-1} r^{-2},
\end{equation*}
hence
\begin{equation}\label{E:la32}
\la \le 32 (1-\ga)^{-2}\ka^{-3} r^{-2}.
\end{equation}

\subsection{}
We are now ready for

\ms
{\bf Proof of Theorem \ref{T:main}.}

The lower bound for $\la$ was established in Lemma \ref{L:lower}.

\ms
We proved the estimate \eqref{E:la32} above under the condition that
there exist $\ga\in (0,1)$, a ball $B_r$ and a regular compact set $F\subset \bB_{3r/2}$
(here the balls $B_r$ and $B_{3r/2}$ have the same center),   such that
$F$ is $\ga$-negligible and its interior includes $\bB_r\setminus \Om$.
(The estimate then holds with $\ka=\ka(\ga,n)$ given by \eqref{E:ka}.)
It follows in particular that $\bB_r\setminus\Om$ is $\ga$-negligible.

Conversely, if $\bB_r\setminus\Om$ is $\ga$-negligible,
then we can approximate it by regular compact sets $F_k$,
$k=1,2,\dots$, such that $\Int F_k\supset \bB_r\setminus \Om$,
$\Int F_k\supset F_{k+1}$, and $\bB_r\setminus \Om$
is the intersection of all $F_k$'s. In this case 
\begin{equation*}
\lim_{k\to\infty}\capa(F_k)=\capa(\bB_r\setminus\Om),
\end{equation*}
due to the continuity property of the capacity. (See e.g. \cite{Mazya8}, Sect. 2.2.1.)
In this case, for any $\eps>0$ the sets $F_k$ will be $(\ga+\eps)$-negligible
for sufficiently large $k$. It follows that the estimate \eqref{E:la32} will hold
if we only know that there exists a  ball $B_r$ such that $\bB_r\setminus\Om$ 
is $\ga$-negligible. Then the estimate still holds  if we replace $r$ by the 
least upper bound of the radii of such balls which is exactly the interior
capacitary radius $r_{\Om,\ga}$. This proves the upper bound in \eqref{E:main}
with
\begin{equation}\label{E:C-ga-n}
C(\ga,n)=32 (1-\ga)^{-2}\ka^{-3},
\end{equation}  
where $\ka$ is defined by \eqref{E:ka}. 
$\square$

\section{Further remarks}

\subsection{Measure instead of capacity}

E.~Lieb \cite{Lieb} used geometric  arguments to establish
a lower bound for $\la(\Om)$ which is similar to 
\eqref{E:lower} but with capacity replaced by the Lebesgue measure.
Such lower bounds can be also deduced from
Theorem \ref{T:main} if we use isoperimetric inequalities
between the capacity and Lebesgue measure:
\begin{equation}\label{E:mes-cap}
\mes F \le A_n (\capa(F))^{n/(n-2)},
\end{equation}
with the equality for balls (see e.g. \cite{Polya} or Sect. 2.2.3, 2.2.4 in \cite{Mazya8}),
so
\begin{equation*}\label{E:An}
A_n=(\mes B_1) \left[\capa(B_1)\right]^{-n/(n-2)}
= 
n^{-1}(n-2)^{-n/(n-2)}\om_n^{-2/(n-2)}.
\end{equation*}
Namely, let us denote for any $\al\in(0,1)$
\begin{equation*}\label{E:r-mes}
r_{\Om,\al}^{(mes)}
=
\sup\{r|\, \exists B_r, \  \mes (B_r\setminus \Om) \le \al\mes B_r\}.
\end{equation*}
Then \eqref{E:mes-cap} implies that
\begin{equation*}\label{E:r-r}
r_{\Om,\al}^{(mes)} \ge r_{\Om,\ga} 
\quad
\text{provided}
\quad
\al = \ga^{n/(n-2)}.
\end{equation*}
Therefore, we obtain
\begin{proposition}\label{P:mes}
For every $\al\in (0,1)$
\begin{equation*}\label{E:mes}
\la(\Om)\ge
c(\ga, n)  \left(r_{\Om,\al}^{(mes)}\right)^{-2},
\quad \text{where} \quad
\ga=\al^{(n-2)/n}.
\end{equation*}
Here $c(\ga,n)$ is given by \eqref{E:c-explicit}.
\end{proposition}
This is exactly Lieb's  inequality (1.2) in \cite{Lieb}, 
though with a different constant.

There are numerous results which give lower bounds for 
$\la(\Om)$. We will mention only a few.  The famous Faber-Krahn  inequality
(\cite{Faber, Krahn, Polya}) gives a lower bound of $\la(\Om)$
in terms of the area of $\Om\subset \R^2$. 
Under miscellaneous topological and geometric restrictions on $\Om$
the interior radius was shown to
provide a lower bound (hence a two-sided estimate) for $\la(\Om)$ in case $n=2$ by
Hayman \cite{Hayman},  Osserman \cite{Osserman-77, Osserman-78, Osserman-79},
Taylor \cite{Taylor}, 
%Cheng \cite{Cheng}, 
Croke \cite{Croke}, Ba\~nuelos and Carroll 
\cite{Banuelos}, and also in case $n\ge 3$
(\cite{Hayman}, \cite{Osserman-79}).

Let $\capa_{\Om}(F)$ denote the capacity of a compact set $F\subset \Om$
with respect to an open set $\Om\subset \R^n$. It is defined similarly to $\capa(F)$ in 
\eqref{E:cap-def}, except the allowed test functions $u$ should be supported in $\Om$.
The following 2-sided estimate for $\la(\Om)$ was established 
in \cite{Mazya-62, Mazya-64}: 
\begin{equation}\label{E:Mazya-Cheeger}
\frac{1}{4}\inf_F \frac{\capa_\Om(F)}{\mes F}
\le 
\la(\Om)
\le
\inf_F \frac{\capa_\Om(F)}{\mes F},
\end{equation}
where the infimum is taken over  all compact sets $F\subset\Om$.
The constant 1/4  in the lower bound is precise.  
Both inequalities  in \eqref{E:Mazya-Cheeger} 
hold on Riemannian manifolds as well.

The lower bound (the first inequality) in \eqref{E:Mazya-Cheeger}
implies the Cheeger inequality \cite{Cheeger} (see Sect. 3 in \cite{Mazya-91} 
for this implication),
which also provides a geometric lower bound for $\la(\Om)$ on manifolds.
(See also Grigor'yan \cite{Grigoryan} for a review and related results.)

\subsection{Bounds for essential spectrum}

Let $\la_\infty(\Om)$ denote the bottom of the essential spectrum of $-\De$
with the Dirichlet boundary conditions in $\Om$.
Then Persson's arguments \cite{Persson} give
\begin{equation*}\label{E:Persson}
\la_\infty(\Om)
= \lim_{R\to+\infty} \la(\Om\setminus \bB_R(0)),
\end{equation*}
 where $B_R(0)$ is the ball with the radius $R$ and the center at the origin.
 Applying two-sided estimates from Theorem \ref{T:main} to 
 $ \la(\Om\setminus \bB_R(0))$, we obtain

\begin{theorem}\label{T:ess}
For any $\ga\in (0,1)$ and any open set $\Om\subset \R^n$,
\begin{equation*}\label{E:ess}
c r_{\Om,\ga,\infty}^{-2} 
\le
\la_\infty(\Om)
\le
C r_{\Om,\ga,\infty}^{-2},
\end{equation*}
 where 
 \begin{equation*}\label{E:r-inf}
 r_{\Om,\ga,\infty}=
 \lim_{R\to\infty} r_{\Om\setminus \bB_R(0),\ga},
 \end{equation*}
and the constants $c=c(\ga,n)$, $C=C(\ga,n)$ are the same as in
Theorem \ref{T:main}.
\end{theorem} 

For small $\ga$ this theorem is due to Maz'ya and Otelbaev (see \cite{Mazya-Otelbaev}
and also Theorem 12.3.1 in \cite{Mazya8}). 

Theorem \ref{T:ess} implies that for any $\ga\in (0,1)$ the condition 
$r_{\Om,\ga,\infty}=0$ is necessary and sufficient for the discreteness
of spectrum of the operator $-\De$ with the Dirichlet boundary
conditions in $\Om$. (This is also a particular case of the main results 
of \cite{Mazya-Shubin}).

\subsection{Bounds for spectra of Schr\"odinger operators}

Theorem \ref{T:main}, Proposition \ref{P:mes} and Theorem \ref{T:ess}
can be extended to Schr\"odinger operators with positive potentials 
(which are even allowed to be  positive measures,
 which are absolutely continuous with respect to 
the Wiener capacity). For small $\ga$ these results can be found in Chapters 10 -- 12
of \cite{Mazya8} with appropriate references. 

For simplicity of formulation we will consider operators $H_V=-\De +V$, $V\ge 0$,
where $V$ is locally integrable. Then 2-sided estimates of the type \eqref{E:main}
can be obtained for the bottom of the spectrum (and essential spectrum)
of $H_V$, if $r_{\Om,\ga}$ is replaced by the quantity 
\begin{equation*}
r_{V,\ga}=
\sup\left\{r|\, \exists B_r, \quad \text{such that}\quad r^{n-2}\ge 
\inf_{F} \int_{B_r\setminus F} V dx\right\},
\end{equation*}
where the infimum is taken over all negligible subsets in $F\subset \bB_r$,
i.e. sets satisfying \eqref{E:ga-negl}.


\begin{thebibliography}{999}

\bibitem{Banuelos} Ba\~nuelos, R., Carroll, T.,
\emph{An improvement of the Osserman constant for the bass note of a drum},
Stochastic Anal. (Ithaca, NY, 1993). Proc. Symp. Pure Math., \textbf{57}, Amer. Math. Soc.,
Providence, RI, 1995, 3--10

\bibitem{Cheeger} Cheeger, J.,
\emph{A  lower bound for the smallest eigenvalue of the Laplacian,}
In: Problems in Analysis, a Symposium in Honor of Salomon Bochner,
Gunning, R.C. (ed.), 
Princeton University Press, Princeton, 1970, 195--199 

\bibitem{Croke} Croke, C.B.,
\emph{The first eigenvalue of the Laplacian for plane domains},
Proc. Amer. Math. Soc., \textbf{81} (1981), 304--305

\bibitem{Faber} Faber, C.,
\emph{Beweis, dass unter allen homogenen Membranen von gleicher Fl\"ache
und gleicher Spannung die kreisf\"ormige den tiefsten Grundton gibt},
Sitzungsber. Bayer. Acad. der Wiss. Math. Phys., Munich 1923, 169--172
 
\bibitem{Grigoryan} Grigor'yan, A.,
\textit{Isoperimetric inequalities and capacities on Riemannian manifolds.}  
The Maz'ya anniversary collection, Vol. 1 (Rostock, 1998).
Oper. Theory Adv. Appl., \textbf{109}, Birkh\"auser, Basel, 1999, 139--153

\bibitem{Hayman} Hayman, W.K.,
\emph{Some bounds for principal frequency,}
Applic. Anal., \textbf{7} (1977/1978), 247--254

\bibitem{Kac} Kac, M.,
\emph{Can one hear the shape of a drum?},
Amer. Math. Monthly, \textbf{73} 1966,  no. 4, part II, 1--23

\bibitem{Krahn} Krahn, E.,
\emph{\"Uber eine von Rayleigh formulierte Minimaleigenschaft des Kreises,}
Math. Ann., \textbf{94} (1925), 97--100

 \bibitem{Lieb} Lieb, E.H.,
\emph{On the lowest eigenvalue of the Laplacian 
for the intersection of two domains,}
Invent. math., \textbf{74} (1983), 441--448

\bibitem{Mazya-62}
Maz'ya, V.G., \textit{The negative spectrum of the higher-dimensional
Schršdinger operator.} (Russian) Dokl. Akad. Nauk SSSR, \textbf{144}, no. 4 (1962), 
721--722. Engl. Transl.: Soviet Math. Dokl., \textbf{3} (1962), 808--810

\bibitem{Mazya-63b}
Maz'ya, V. G., \textit{On the boundary regularity of solutions of elliptic
equations and of a conformal mapping.} (Russian) Dokl. Akad. Nauk SSSR,
\textbf{152} (1963), 1297--1300 

\bibitem{Mazya-64}  Maz'ya, V.G.,  \textit{On the theory of the multidimensional
Schr\"odinger operator.} (Russian) Izv. Akad. Nauk SSSR Ser. Mat., \textbf{28} (1964),
1145--1172

\bibitem{Mazya-74}
Maz'ya, V. G.,  \textit{The connection between two forms of capacity.} (Russian)
Vestnik Leningrad. Univ. Mat. Mech. Astronom., \textbf{7}:2 (1974), 33--40


\bibitem{Mazya8} Maz'ya, V.G.,
\emph{Sobolev spaces},
Springer Verlag, Berlin, 1985

\bibitem{Mazya-91} Maz'ya, V.G.,
\emph{Classes of domains, measures and capacities in the theory
of differentiable functions}. Analysis III. Spaces of Differentiable Functions.
Encyclopaedia of Math. Sciences, vol. 26, Springer-Verlag, 1991, 141--211

\bibitem{Mazya-Otelbaev} Maz'ya, V.G., Otelbaev, M., 
\emph{Imbedding theorems and the spectrum of a certain pseudodifferential operator.} 
(Russian)  Sibirsk. Mat. Z., \textbf{18} (1977), no. 5, 1073--1087

\bibitem{Mazya-Shubin}  Maz'ya, V.G., Shubin, M.A.,
\emph{Discreteness of spectrum and positivity criteria for Schr\"odinger operators.}
%Preprint math.SP/0305278. 
Annals of Math., \textbf{162} (2005), 919--942

\bibitem{Osserman-77} Osserman, R.,
\emph{A note on Hayman's theorem on the bass note of a drum},
Comment. Math. Helv., \textbf{52} (1977),  545--555

\bibitem{Osserman-78} Osserman, R.,
\emph{The isoperimetric Inequality}, 
Bull. Amer. Math. Soc., \textbf{84} (1978), 1182--1238

\bibitem{Osserman-79} Osserman, R.,
\emph{Bonnesen-style isoperimetric inequalities},
Amer. Math. Monthly, \textbf{86} (1979), 1--29


\bibitem{Persson} Persson, A.,
\emph{Bounds for the discrete part of the spectrum of a semi-bounded Schr\"odinger operator},
Math. Scand., \textbf{8} (1960), 143--153

\bibitem{Polya} P\'olya, G., Szeg\"o, G.,
\emph{Isoperimetric inequalities in mathematical physics}, 
Princeton University Press, Princeton,  1951.

\bibitem{Rogers} Rogers, C.A.,
\emph{Packing and covering},
Cambridge University Press, 1964

\bibitem{Taylor} Taylor, M.,
\emph{Estimate on the fundamental frequency of a drum},
Duke Math. J., \textbf{46} (1979), 447--453 

\bibitem{Wermer}   Wermer, J.,  \emph{Potential theory}, 
Lecture Notes in Math., \textbf{408}, Springer-Verlag, 1974 

\end{thebibliography}
\end{document}